# A non-linear Renewal Theorem with stationary and slowly changing perturbations

## Dong-Yun Kim[1] and Michael Woodroofe[2]

*Michigan State University and the University of Michigan*

**Abstract:** Non-linear renewal theory is extended to include random walks perturbed by both a slowly changing sequence and a stationary one. Main results include a version of the Key Renewal Theorem, a derivation of the limiting distribution of the excess over a boundary, and an expansion for the expected first passage time. The formulation is motivated by problems in sequential analysis with staggered entry, where subjects enter a study at random times.

## 1. Introduction

Non-linear renewal theory concerns sequences of the form $Z_n^o = S_n + \zeta_n$, where $S_n$ is a random walk with a finite positive drift, $\mu$ say, and $\zeta_n$ is a sequence that changes slowly, as in (5) and (6) below. The main results are extensions of the Renewal Theorem and selected corollaries from the random walk $S_n$ to the perturbed random walk $Z_n^o$. For example, letting $t_a^o$ denote the first passage time $t_a^o = \inf\{n \geq 1 : Z_n^o > a\}$, it is shown that $E(t_a^o) = (a + \rho - \eta)/\mu + o(1)$ under natural conditions, where $\rho$ is the mean of the asymptotic distribution of the excess $Z_{t_a^o} - a$ and $\eta = \lim_{n \to \infty} E(\zeta_n)$. The formulation is due to Lai and Siegmund [5] and [6]. The early development of non-linear renewal theory is described in Woodroofe [8], and more general results have since been obtained by Zhang [9]. The formulation is motivated by problems from sequential analysis, where the stopping times of many sequential procedures may be written in the form $t_a^o$. Applications to sequential analysis are described in [8] and in Siegmund's book [7].

Sequential problems with staggered entry (where subjects arrive according to a Poisson process, say) lead to random walks perturbed by both a slowly changing term and a stationary one, which is not slowly changing. Examples are described below. In such cases, there is interest in processes of the form

$$(1) \qquad Z_n = S_n + \xi_n + \zeta_n,$$

where $S_n$ and $\zeta_n$ are as above and $\xi_n$ is a stationary sequence and in the associated first passage times,

$$(2) \qquad t_a = \inf\{n \geq n_0 : Z_n > a\},$$

---

[1]Department of Probability and Statistics, Michigan State University, A413 Wells Hall, East Lansing, MI 48824-1027, USA, e-mail: `dhkim@ilstu.edu`

[2]Department of Statistics, University of Michigan, 439 West Hall, 1085 South University, Ann Arbor, MI 48109-1107, e-mail: `michaelw@umich.edu`

*AMS 2000 subject classifications:* 60K05.

*Keywords and phrases:* excess over the boundary, first passage times, fixed width confidence intervals, repeated likelihood ratio tests, staggered entry.

where $n_0$ is a fixed (throughout) positive integer. Kim and Woodroofe [4] considered a simple case of (1) in which $\zeta_n \equiv 0$, obtaining versions of the Renewal Theorem and its corollaries in that context. They also applied their results to the problem of testing a simple hypothesis about the mean of an exponential distribution, but applications were severely limited by the absence of a slowly changing term. The purpose of this paper is to supply a non-linear renewal theory in the context of (1).

For example, consider a model in which patients arrive according to a Poisson process, with rate $\lambda > 0$ say, are treated, and then live for an exponentially distributed residual life $L_k$ with unknown failure rate $\theta > 0$, (or are cured after an exponentially distributed time with unknown success rate). Letting $\tau_{k-1}$ denote the arrival time of the $k^{\text{th}}$ patient with $\tau_0 = 0$, the data available at time $\tau_n$ are $L_{n,k} = \min[L_k, \tau_n - \tau_{k-1}]$ and $\delta_{n,k} = \mathbf{1}\{L_k \leq \tau_n - \tau_{k-1}\}$, $k = 1, \ldots, n$, and the log-likelihood function is

$$\ell_n(\theta) = K_n \log(\theta) - \theta T_n^*,$$

where $K_n = \delta_{n,1} + \cdots + \delta_{n,n}$ and $T_n^* = L_{n,1} + \cdots + L_{n,n}$. Here $K_n$ counts the number of deaths (or cures), and $T_n^*$ is called the *total time on test statistic*. To see how staggered entry leads to processes of the form (1), write

$$T_n^* = \sum_{k=1}^n L_k - \sum_{k=1}^n [L_{n-k+1} - (\tau_n - \tau_{n-k})]_+.$$

The first term on the right is a random walk, and the second may be written as

$$\sum_{k=1}^n [L_{n-k+1} - (\sum_{j=0}^{k-1} \eta_{n-j})]_+ = \sum_{k=1}^\infty [L_{n-k+1} - (\sum_{j=0}^{k-1} \eta_{n-j})]_+ + o(1) \ w.p.1,$$

where $\eta_k = \tau_k - \tau_{k-1}$. If the (sequences of) arrival times and lifetimes are independent, then the last sum is a stationary sequence and the $o(1)$ is slowly changing. So, $T_n^*$ is of the form (1). Combining these observations with Taylor series expansions, leads to processes of the form (1) for many other statistics. Specific examples are considered in Section 4.

The main results are stated and proved in Section 3. Section 2 contains some preliminary lemmas.

## 2. Preliminaries

In Sections 2 and 3, $\ldots W_{-1}, W_0, W_1, W_2, \ldots$ denote i.i.d. random elements with values in a Polish space, $\mathcal{W}$ say, $X_k = \varphi(W_k)$ and $Y_k = \psi(W_k)$, where $\varphi : \mathcal{W} \to \mathbb{R}$ and $\psi : \mathcal{W} \to \mathbb{R}^d$ are Borel measurable functions, $S_n = X_1 + \cdots + X_n$, and $T_n = Y_1 + \cdots + Y_n$. Suppose throughout that $X_1$ has a non-arithmetic distribution with a finite, positive mean $\mu$ and that $E(Y_k) = 0$. Next, let $\xi_n$ be a sequence of the form

(3) $$\xi_n = \xi(W_n, W_{n-1}, \ldots),$$

where $\xi$ is measurable on $\mathcal{W}^{\mathbb{N}}$ and $\mathbb{N} = \{0, 1, 2, \ldots\}$. Let $\zeta_n = \zeta_n(W_1, \ldots, W_n)$ be random variables of the form

(4) $$\zeta_n = \zeta_n' + \zeta_n''$$

where

$$\zeta'_n = \frac{1}{n} T'_n Q T_n, \tag{5}$$

$Q$ is a symmetric $d \times d$ matrix, and

$$\lim_{n \to \infty} n^\kappa \sum_{n \leq k \leq n+n^{\frac{2}{3}}} P(|\zeta''_k| \geq \epsilon) = 0 \tag{6}$$

for each $\epsilon > 0$ for appropriate $\kappa \geq 0$. Finally, suppose (throughout) that

$$E|X_1|^{p_x} + E|\xi_1|^{p_\xi} + E\|Y_1\|^{2p_y} + \sup_{n \geq n_0} E|\zeta_n|^{p_\zeta} < \infty, \tag{7}$$

for appropriate $p_x$, $p_y$, $p_\xi$, and $p_\zeta \geq 2$. Then

$$\zeta_n \Rightarrow \zeta = \sum_{i=1}^d \lambda_i \chi^2_{1,i} \tag{8}$$

where $\chi^2_{1,1}, \ldots, \chi^2_{1,d}$ are independent chi-square random variables and $\lambda_1, \ldots, \lambda_d$ are constants. Let $L$ denote the distribution function of $\zeta$. Then $L$ is continuous unless $Y'_1 Q Y_1 = 0$ w.p.1. To avoid trivialities we suppose that $P[Y'_1 Q Y_1 \neq 0] > 0$ for the remainder of the paper.

The conditions imposed above with $\kappa = 0$ and $p_x = p_y = p_\xi = p_\zeta = 2$ are *Standing Assumptions* and are not repeated in the statements of results. Additional moment conditions are imposed by requiring (7) for higher values of $p_x$, $p_y$, $p_\xi$, and $p_\zeta$, and (6) for selected $\kappa > 0$.

With $Z_n = S_n + \xi_n + \zeta_n$, as in (1), and $1/3 < q < 1/2$, let

$$m = m_a = \left\lfloor \left( \frac{1 - a^{-q}}{\mu} \right) a \right\rfloor, \tag{9}$$

$$M = M_a = \left\lfloor \left( \frac{1 + a^{-q}}{\mu} \right) a \right\rfloor, \tag{10}$$

$$\Delta_0(q, a) = \sum_{n=0}^m P(Z_n > a), \tag{11}$$

and

$$\Delta_1(q, a, b) = \sum_{n=M+1}^\infty P(Z_n \leq a + b), \tag{12}$$

where $\lfloor x \rfloor$ is the greatest integer that is less than or equal to $x$. Below, there is special interest in $q$'s that are slightly larger than $1/3$.

Lemmas 1 and 2 are extensions and improvements of Lemmas 1 and 4 in [4]. The proof of Lemma 2 is virtually identical to the proof of Lemma 4 of [4]. It is only necessary to replace the $\mathbf{W}_n$ of [4] by $\mathbf{W}_n^*$. The proof of Lemma 1 is outlined in the Appendix.

**Lemma 1.** *If $p_x \geq (2-q)/(1-q)$, then*

$$\text{(13)} \qquad \lim_{a\to\infty}\left\{\Delta_0(q,a) + \Delta_1(q,a,\frac{1}{2}a^{1-q}) + \sum_{n=M}^{\infty} P[t_a > n]\right\} = 0.$$

*If $p_x \geq 2/(1-q)$, $p_\xi \wedge p_\zeta > (2-q)/(1-q)$ then*

$$\text{(14)} \qquad \lim_{a\to\infty}\left\{aP[t_a \leq m] + aP[t_a > M]\right\} = 0.$$

For Lemmas 2 and 4, let $\mathcal{B}$ denote the Borel sets of $\mathcal{W}$, $\mathcal{B}^k = \mathcal{B} \otimes \cdots \otimes \mathcal{B}$ ($k$ factors), $\mathcal{B}^{\mathbb{N}}$ the product sigma-algebra in $\mathcal{W}^{\mathbb{N}}$, and $\mathcal{B}^{\mathbb{Z}}$ the product sigma-algebra in $\mathcal{W}^{\mathbb{Z}}$, where $\mathbb{Z} = \{\ldots -1,0,1,\ldots\}$. Further, let

$$\mathbf{W}_n = (\ldots W_{n-1}, W_n),$$
$$\mathbf{W}_n^* = (\ldots W_{n-1}, W_n, W_{n+1}, \ldots),$$
$$T_{m,n} = \sum_{k=n-m+1}^{n} Y_k,$$
$$\tilde{\zeta}_{m,n} = \frac{1}{m}T'_{m,n}QT_{m,n},$$

and

$$\tilde{Z}_{m,n} = S_n + \xi_n + \tilde{\zeta}_{m,n}.$$

Thus, the $j^{\text{th}}$ component of $\mathbf{W}_n^*$ is $W_{n+j}$. Observe that $\xi_n + \tilde{\zeta}_{m,n}$ is a stationary process (in $n$) for fixed $a$. So, $\tilde{Z}_{m,n}$ is a random walk perturbed only by a stationary sequence (that depends on $a$).

**Lemma 2.** *Given $0 < \alpha < 1$, let $\beta = (1-\alpha^2)/\alpha$. Then there are constants $K_1 = K_{1,\alpha}$ and $K_2 = K_{2,\alpha}$ for which*

$$\sum_{n=0}^{\infty} P[\mathbf{W}_n^* \in B,\ a < S_n \leq a+b] \leq K_1(1+b)P[\mathbf{W}_0^* \in B]^{\alpha}$$

*and*

$$\sum_{n=0}^{\infty} P[\mathbf{W}_n^* \in B,\ a < \tilde{Z}_{m,n} \leq a+b] \leq K_2(1+b)\left[\int_{\mathbf{W}_0^* \in B}(1 + |\xi_0 + \tilde{\zeta}_{m,0}|)^{\beta} dP\right]^{\alpha}$$

*for all $B \in \mathcal{B}^{\mathbb{Z}}$ and $0 \leq a, b < \infty$.*

**Lemma 3.** *Let $m$, $M$ be as in (9) and (10) and $\epsilon > 0$. Then*

$$\text{(15)} \qquad \lim_{a\to\infty}\sum_{n=m+1}^{M} P[|\zeta_n - \tilde{\zeta}_{m,n}| \geq \epsilon] = 0;$$

*and if $\kappa = 1$ in (6) and $p_y > 1/q$ in (7), then*

$$\text{(16)} \qquad \lim_{a\to\infty} aP[\max_{m \leq n \leq M}|\zeta_n - \tilde{\zeta}_{m,n}| \geq \epsilon] = 0.$$

*Proof.* By (6), it suffices to establish (15) with $\zeta_n$ replaced by $\zeta'_n$. Towards this end observe that

$$\zeta'_n - \tilde{\zeta}_{m,n} = \left(\frac{1}{n} - \frac{1}{m}\right) T'_n Q T_n + \frac{1}{m}(T'_n Q T_n - T'_{m,n} Q T_{m,n})$$

and

$$P(|\zeta'_n - \tilde{\zeta}_{m,n}| \geq \epsilon) \leq I_n + II_n$$

where

$$I_n = P\left(\left|\frac{1}{n} - \frac{1}{m}\right| |T'_n Q T_n| \geq \frac{\epsilon}{2}\right)$$

and

$$II_n = P\left(\frac{1}{m}|T'_n Q T_n - T'_{m,n} Q T_{m,n}| \geq \frac{\epsilon}{2}\right).$$

Since $Y_k$ have finite fourth moments, there is a constant $C$ for which $E[(T'_n Q T_n)^2] \leq Cn^2$ and, therefore

$$I_n \leq P\left(|T'_n Q T_n| \geq \frac{\epsilon m^2}{2(M-m)}\right) \leq \frac{4(M-m)^2}{\epsilon^2 m^4} \cdot Cn^2 \leq \frac{8C(M-m)^2}{\epsilon^2 m^2}$$

for all $m \leq n \leq M$ and all sufficiently large $a$. For $II_n$, first observe that $T_n = T_{m,n} + T_{n-m}$ where $T_{m,n}$ and $T_{n-m}$ are independent, and that $T_{m,n}$ has the same distribution as $T_m$. So,

$$II_n \leq P[2||T_{m,n}|| \, ||Q|| \, ||T_{n-m}|| \geq \frac{m\epsilon}{4}] + P[||Q|| \, ||T_{n-m}||^2 \geq \frac{m\epsilon}{4}]$$
$$\leq \frac{C'm^2(n-m)^2}{\epsilon^4 m^4} + \frac{C'(n-m)^2}{\epsilon^2 m^2}$$
$$\leq \frac{2C'}{\epsilon^4} \cdot \frac{(M-m)^2}{m^2}.$$

for $0 < \epsilon < 1$, where $C'$ is a constant depending on the 4th moment of $||Y_1||$ and $||Q||$. Since $q > 1/3$,

$$\sum_{n=m+1}^{M}[I_n + II_n] \leq \frac{8C + 2C'}{\epsilon^4} \frac{(M-m)^3}{m^2} \to 0,$$

establishing (15).

If (6) holds with $\kappa = 1$, then $P[\max_{m \leq n \leq M} |\zeta''_n| \geq \epsilon] = o(1/a)$ for each $\epsilon > 0$ by Boole's Inequality, and it suffices to establish (16) with $\zeta_n$ replaced by $\zeta'_n$. For this

$$P[\max_{m \leq n \leq M} |\zeta'_n - \tilde{\zeta}_{m,n}| \geq \epsilon] \leq I^*_a + II^*_a$$

where

$$I^*_a = P[\max_{m \leq n \leq M} |T'_n Q T_n| \geq \frac{m^2 \epsilon}{2(M-m)}]$$

and

$$II^*_a = P[\max_{m \leq n \leq M} 2|T'_{m,n} Q T_{n-m}| \geq \frac{m\epsilon}{4}] + P[\max_{m \leq n \leq M} |T'_{n-m} Q T_{n-m}| \geq \frac{m\epsilon}{4}].$$

By the Submartingale Inequality, there is a constant $C''$ for which

$$I_a^* \leq C''\big[\frac{2(M-m)}{m^2\epsilon}\big]^{p_y} M^{p_y},$$

which is $o(1/a)$ if $p_y > 1/q$. For $II_a^*$, observe that $T_{m,n} = \sum_{k=n-m+1}^{M-m} Y_k + \sum_{k=M-m+1}^{n} Y_k$. So, using the Submartingale Inequality again and possibly enlarging $C''$ and

$$II_a^* \leq C''\big[\frac{4}{m\epsilon}\big]^{2p_y} M^{p_y}(M-m)^{p_y} + C''\big[\frac{4}{m\epsilon}\big]^{p_y}(M-m)^{p_y},$$

which is again $o(1/a)$ if $p_y > 1/q$.  □

Let $\mathcal{C}_k$ the collection of cylinder subsets in $\mathcal{W}^{\mathbb{N}}$ with base in $\mathcal{W}^k$, and write members of $\mathcal{C}_k$ as $\mathcal{W}^{\mathbb{N}} \times C$, where $C \in \mathcal{W}^k$. Also let $\mathcal{C} = \cup_{k=1}^{\infty}\mathcal{C}_k$ be the class of cylinder sets.

**Lemma 4.** *If $p_x > (2-q)/(1-q)$, then for each $0 < b, z < \infty$,*

$$\lim_{a\to\infty} |\sum_{n=m+1}^{\infty} P[\mathbf{W}_n \in B, \zeta_n \leq z, a < S_n \leq a+b]$$
(17)
$$- \frac{b}{\mu}P[\mathbf{W}_0 \in B, \tilde{\zeta}_{m,n} \leq z]| = 0$$

*uniformly with respect to $B \in \mathcal{C}_m$; and*

(18) $$\lim_{a\to\infty} \sum_{n=1}^{\infty} P[\mathbf{W}_n \in B, \zeta_n \leq z, a < S_n \leq a+b] = \frac{b}{\mu}P[\mathbf{W}_0 \in B]L(z)$$

*for each $B \in \mathcal{B}^{\mathbb{N}}$.*

*Proof.* First consider (17) when $z = \infty$ and let $V$ denote the renewal function for $S_1, S_2, \ldots$. If $B \in \mathcal{C}_m$, so that $B = \mathcal{W}^{\mathbb{N}} \times C$ for some $C \in \mathcal{B}^m$, then

$$\sum_{n=m+1}^{\infty} P(\mathbf{W}_n \in B, \ a < S_n \leq a+b) = \int_{\mathbf{W}_m \in B} V(a - S_m, a - S_m + b)dP,$$

by independence and symmetry, as in Lemma 2 of [4]. As $a \to \infty$, $a - S_m = (a - m\mu) + (m\mu - S_m) \to \infty$ in probability since $q < 1/2$ and, therefore, $V(a - S_m, a - S_m + b) \to b/\mu$ in probability, by the Renewal Theorem for random walks, [3], Ch. 12. It follows that

$$|\int_{\mathbf{W}_m \in B} V(a - S_m, a - S_m + b)dP - \frac{b}{\mu}P[\mathbf{W}_m \in B]|$$
$$\leq \int |V(a - S_m, a - S_m + b) - \frac{b}{\mu}|dP \to 0$$

uniformly with respect to $B \in \mathcal{C}_m$. That is, (17) holds when $z = \infty$, since $P[\mathbf{W}_m \in B] = P[\mathbf{W}_0 \in B]$.

Now consider the case $z < \infty$. Using Lemmas 1 and 3,

$$\sum_{n=m+1}^{\infty} P[\mathbf{W}_n \in B, \tilde{\zeta}_{m,n} \leq z - \epsilon, \ a < S_n \leq a + b] - o(1)$$

$$\leq \sum_{n=m+1}^{\infty} P[\mathbf{W}_n \in B, \ \zeta_n \leq z, \ a < S_n \leq a + b]$$

$$\leq \sum_{n=m+1}^{\infty} P[\mathbf{W}_n \in B, \ \tilde{\zeta}_{m,n} \leq z + \epsilon, \ a < S_n \leq a + b] + o(1).$$

for any $\epsilon > 0$, where the $o(1)$ terms are independent of $B \in \mathcal{B}^{\mathbb{N}}$. Let

$$B_a^{\pm} = \{\mathbf{w} \in B : \tilde{\zeta}_{m,0}(\mathbf{w}) \leq z \pm \epsilon\}.$$

Then the summands on the first and last lines are $P[\mathbf{W}_n \in B_a^{\pm}, \ a < S_n \leq a + b]$, and $P[\mathbf{W}_0 \in B_a^{\pm}] = P[\mathbf{W}_0 \in B, \ \tilde{\zeta}_{m,0} \leq z \pm \epsilon]$. If $B \in \mathcal{C}_m$, then also $B_a^{\pm} \in \mathcal{C}_m$ and, therefore,

$$\lim_{a \to \infty} \left\{ \sum_{n=m+1}^{\infty} P[\mathbf{W}_n \in B_a^{\pm}, \ a < S_n \leq a + b] - \frac{b}{\mu} P[\mathbf{W}_0 \in B, \tilde{\zeta}_{m,0} \leq z \pm \epsilon] \right\} = 0,$$

uniformly with respect to $B \in \mathcal{C}_m$. It follows that

$$\limsup_{a \to \infty} |\sum_{n=m+1}^{\infty} P[\mathbf{W}_n \in B, \ \tilde{\zeta}_{m,n} \leq z, a < S_n \leq a + b] - \frac{b}{\mu} P[\mathbf{W}_0 \in B, \ \tilde{\zeta}_{m,0} \leq z]|$$

$$\leq \limsup_{a \to \infty} P[z - \epsilon \leq \tilde{\zeta}_{m,0} \leq z + \epsilon] = L(z + \epsilon) - L(z - \epsilon)$$

uniformly with respect to $B \in \mathcal{C}_m$ and this establishes (17), since $\epsilon > 0$ was arbitrary.

Now consider (18). For any $B \in \mathcal{B}^{\mathbb{N}}$, and $\epsilon > 0$, there is a cylinder set $B_0 \in \mathcal{C}$ for which $P[\mathbf{W}_0 \in B \Delta B_0] \leq \epsilon$, where $\Delta$ denotes symmetric difference, in which case

$$|\sum_{n=1}^{\infty} P[\mathbf{W}_n \in B, \ \tilde{\zeta}_{m,0} \leq z, \ a < S_n \leq a + b]$$

$$- P[\mathbf{W}_n \in B_0, \ \tilde{\zeta}_{m,0} \leq z, \ a < S_n \leq a + b]| \leq K_1(1 + b)\sqrt{\epsilon},$$

by Lemma 2. Combining this inequality with (17),

$$\limsup_{a \to \infty} |\sum_{n=1}^{\infty} P[\mathbf{W}_n \in B, \ \tilde{\zeta}_{m,n} \leq z, \ a < S_n \leq a + b]$$

$$- \frac{b}{\mu} P[\mathbf{W}_0 \in B_0, \ \tilde{\zeta}_{m,0} \leq z]| \leq K_1(1 + b)\sqrt{\epsilon}.$$

If $B \in \mathcal{C}_k$, write $T_{m,0} = (Y_{-m+1} + \cdots + Y_{-k-1}) + T_{k+1,0}$. Here $B_0$ and the first term are independent and the second is bounded as $m \to \infty$. It follows easily that

$$\lim_{a \to \infty} P[\mathbf{W}_0 \in B_0, \ \tilde{\zeta}_{m,0} \leq z] = P[\mathbf{W}_0 \in B_0]L(z),$$

from which the Theorem follows by letting $\epsilon \to 0$. □

## 3. Main results

There are four main results: a renewal theorem, uniform integrability, the limiting distribution of the excess over the boundary, and properties of the first passage time $t_a$. For technical reasons, they are presented in that order. Recall that $1/3 < q < 1/2$.

**Theorem 1.** *If $p_x > (2-q)/(1-q)$, then*

$$(19) \quad \lim_{a \to \infty} \sum_{n=1}^{\infty} P[\mathbf{W}_n \in B, \zeta_n \leq y, a < Z_n \leq a+b] = \frac{b}{\mu} P[\mathbf{W}_0 \in B] L(y)$$

*for each $B \in \mathcal{B}^{\mathbb{N}}$ and $0 \leq y < \infty$.*

*Proof.* It suffices to establish the theorem for $0 < b \leq \mu$. Let

$$J_a(x, y, z) = \sum_{n=1}^{\infty} P[\mathbf{W}_n \in B, \xi_n \leq x, \zeta_n \leq y, a < S_n \leq a+z]$$

and

$$K_a(x, y, z) = \sum_{n=1}^{\infty} P[\mathbf{W}_n \in B, \xi_n \leq x, \zeta_n \leq y, a < Z_n \leq a+z]$$

for $0 < a, b < \infty$ and $y, z \in \mathbb{R}$. Then

$$(20) \quad \lim_{a \to \infty} J_a(x, y, z) = \frac{z}{\mu} P[\mathbf{W}_0 \in B, \xi_0 \leq x] L(y) := J_\infty(x, y, z)$$

by Lemma 4 with the $B$ in Lemma 4 replaced by $B \cap \{\xi_0 \leq x\}$. Next, let $\Gamma(x, y, z) = (x, y, x+y+z)$, where addition is understood modulo $\mu$. Then $\Gamma$ is continuous a.e. $(J_\infty)$ and $K_a = J_a \circ \Gamma^{-1}$, and $J_\infty \circ \Gamma^{-1} = J_\infty$. Theorem 1 now follows from the Continuous Mapping Theorem. □

The uniform integrability that is needed will be deduced as corollaries to Theorem 2 below. In its statement, $A_a$ denotes the event

$$A_a = \{m < t \leq M, \max_{m \leq n \leq M} |\zeta_n - \tilde{\zeta}_{m,n}| \leq 1\}.$$

In the proofs $U_n = Z_n - Z_{n-1} = X_n + (\xi_n - \xi_{n-1}) + (\zeta_n - \zeta_{n-1})$ and $\tilde{U}_{m,n} = \tilde{Z}_{m,n} - \tilde{Z}_{m,n-1} = X_n + (\xi_n - \xi_{n-1}) + (\tilde{\zeta}_{m,n} - \tilde{\zeta}_{m,n-1})$.

**Theorem 2.** *Suppose $p_x > (2-q)/(1-q)$. If $0 < \alpha < 1$, there is a constant $K_3 = K_{3,\alpha}$ for which*

$$P[A_a, \mathbf{W}_{t_a}^* \in B] \leq K_3 \Big[ \int_{\mathbf{W}_0^* \in B} (1 + |\xi_0 + \tilde{\zeta}_{m,0}|)^\beta (1 + \tilde{U}_{m,0})^{1+\beta} dP \Big]^\alpha,$$

*for all Borel sets $B \in \mathcal{B}^{\mathbb{N}}$, where $\beta = (1-\alpha^2)/\alpha$.*

*Proof.* Since $t_a = n$ implies $Z_{n-1} \leq a < Z_n$,

$$P[A_a, \mathbf{W}_{t_a}^* \in B] \leq \sum_{n=m+1}^{M} P[\mathbf{W}_n^* \in B, \tilde{Z}_{m,n-1} \leq a+1, \tilde{Z}_{m,n} > a-1]$$

$$\leq \sum_{n=m+1}^{M} \sum_{k \leq \lfloor a \rfloor + 2} P[\mathbf{W}_n^* \in B, k-1 < \tilde{Z}_{m,n-1} \leq k, \tilde{U}_{m,n} > a-k-1]$$

$$\leq 2K_2 \sum_{k \leq \lfloor a \rfloor + 2} \Big[ \int_{\mathbf{W_0}^* \in B, \tilde{U}_{m,0} > a-k-1} (1 + |\xi_0 + \tilde{\zeta}_{m,0}|)^\beta dP \Big]^\alpha,$$

using Lemma 2 in the last step. Letting $j = \lfloor a \rfloor - k - 1$, the latter sum is at most

$$8K_2\Big[\int_{\mathbf{W_{o^*}}\in B}(1+|\xi_0+\tilde{\zeta}_{m,0}|)^\beta dP\Big]^\alpha$$
$$+ 2K_2\sum_{j=0}^{\infty}\Big[\int_{\mathbf{W_{o^*}}\in B,\tilde{U}_{m,0}>j}(1+|\xi_0+\tilde{\zeta}_{m,0}|)^\beta dP\Big]^\alpha$$

of which the sum is at most

$$\sum_{j=0}^{\infty}\Big(\frac{1}{1+j}\Big)^{\alpha\beta}\Big[\int_{\mathbf{W_{o^*}}\in B,\tilde{U}_{m,0}>j}(1+|\xi_0+\tilde{\zeta}_{m,0}|)^\beta(1+\tilde{U}_{m,0})^\beta dP\Big]^\alpha$$
$$\leq \Big[\sum_{j=0}^{\infty}\Big(\frac{1}{1+j}\Big)^{\frac{\alpha\beta}{1-\alpha}}\Big]^{1-\alpha}\cdot\Big[\sum_{j=0}^{\infty}\int_{\mathbf{W_{o^*}}\in B,\tilde{U}_{m,0}>j}(1+|\xi_0+\tilde{\zeta}_{m,0}|)^\beta(1+\tilde{U}_{m,0})^\beta dP\Big]^\alpha$$
$$\leq \Big[\sum_{j=0}^{\infty}\Big(\frac{1}{1+j}\Big)^{1+\alpha}\Big]^{1-\alpha}\cdot\Big[\int_{\mathbf{W_{o^*}}\in B}(1+|\xi_0+\tilde{\zeta}_{m,0}|)^\beta(1+\tilde{U}_{m,0})^{1+\beta} dP\Big]^\alpha$$

The theorem follows easily. □

The uniform integrability over $A_a$ of $R_a$, $\xi_{t_a} + \zeta_{t_a}$ is easy consequence.

**Corollary 1** (Under the conditions of Theorem 2). $R_a \mathbf{1}_{A_a}$ *are uniformly integrable.*

*Proof.* If $A_a$ occurs, then $R_a \leq U_{t_a} \leq \tilde{U}_{m,t_a} + 2$, and

$$P[A_a, \tilde{U}_{m,t_a} > r] \leq K_3\Big[\int_{\tilde{U}_{m.0}>r}(1+|\xi_0+\tilde{\zeta}_{m,0}|)^\beta(1+\tilde{U}_{m,0}^+)^{\beta+1}dP\Big]^\alpha$$
$$\leq \frac{K_3}{r^{\alpha\gamma}}\Big[\int_{\tilde{U}_{m.0}>r}(1+|\xi_0+\tilde{\zeta}_{m,0}|)^\beta(1+\tilde{U}_{m,0}^+)^{\beta+\gamma+1}dP\Big]^\alpha$$
$$\leq \frac{K_3}{r^{\alpha\gamma}}\Big[\int_{\tilde{U}_{m.0}>r}(1+|\xi_0+\tilde{\zeta}_{m,0}|)^{2\beta+\gamma+1}dP\Big]^{\frac{\alpha\beta}{2\beta+\gamma+1}}$$
$$\cdot\Big[\int_{\tilde{U}_{m.0}>r}(1+\tilde{U}_{m,0}^+)^{2\beta+\gamma+1}dP\Big]^{\alpha\frac{\beta+\gamma+1}{2\beta+\gamma+1}}$$

for any $0 < \alpha < 1$ and $\gamma \geq 0$, where $\beta = (1-\alpha^2)/\alpha$ and $K_3$ are as in Theorem 2. The conditions of the theorem require $p_x > 5/2$. So, there are $\alpha < 1$ and $\gamma > 1/\alpha$ for which $2\beta + \gamma + 1 \leq p_x$. Then the last two integrals are bounded in $a$ and $r$, and the uniform integrability of $\tilde{U}_{m,t_a}$ over $A_a$ follows. □

**Corollary 2.** $|\xi_{t_a} + \zeta_{t_a}|\mathbf{1}_{A_a}$ *are uniformly integrable.*

*Proof.* Similar to that of Corollary 1. □

The next theorem uses the following easily verified lemma: If $V_a$ and $V$ are random vectors for which (the distributions of) $V_a$ are tight and

(21) $$\limsup_{a\to\infty} P[V_a \in K] \leq P[V \in K]$$

for all compact rectangles $K$, then $V_a$ converges in distribution to $V$

**Theorem 3.** *If $p_x > (2-q)/(1-q)$ then for all $0 < r < s < \infty$,*

$$\lim_{a \to \infty} P[r < R_a \leq s, \xi_{t_a} \leq y, \zeta_{t_a} \leq z] = \frac{1}{\mu} \int_r^s P[\inf_{j \leq -1} Z_j^* \geq u, \xi_0 \leq y] du \times L(z),$$

*where $Z_j^* = X_{j+1} + \cdots + X_0 + \xi_0 - \xi_j$.*

*Proof.* If $I$ and $J$ are intervals, then

$$P[r \leq R_a \leq s, \xi_{t_a} \in I, \zeta_{t_a} \in J]$$
$$= \sum_{n=1}^{\infty} P[t_a \geq n, \xi_n \in I, \zeta_n \in J, a + r \leq Z_n \leq a + s]$$

Next, for any $\epsilon > 0$ and any integer $k$,

$$P[t_a \geq n, \ \xi_n \in I, \zeta_n \in J, \ a + r \leq Z_n \leq a + s]$$
$$\leq P[\min_{j \leq k}(Z_n^1 - Z_{n-j}^1) \geq r - \epsilon, \xi_n \in I, \zeta_n \in J, a + r \leq Z_n \leq a + s]$$
$$+ P[\max_{j \leq k} |\zeta_n - \zeta_{n-j}| > \epsilon],$$

where $Z_n^1 = S_n + \xi_n$; and $\min_{j \leq k}(Z_n^1 - Z_{n-j}^1) \geq r - \epsilon$ iff $\mathbf{W}_n \in B$, where

$$B = \{\mathbf{w}_0 : x_{-j+1} + \cdots + x_0 + \xi(\mathbf{w}_0) - \xi(\mathbf{w}_{-j+1}) \geq r - \epsilon, \text{ for } j = 1, \cdots, k\}.$$

So,

$$P[r \leq R_a \leq s, \xi_{t_a} \in I, \zeta_{t_a} \in J]$$
$$\leq \sum_{n=1}^{\infty} P[\mathbf{W}_n \in B, \xi_n \in I, \zeta_n \in J, a + r \leq Z_n \leq a + s] + o(1),$$

and

$$\limsup_{a \to \infty} P[r \leq R_a \leq s, \xi_{t_a} \in I, \zeta_{t_a} \in J] \leq (\frac{s-r}{\mu}) P[\mathbf{W}_0 \in B, \xi_0 \in I] L\{J\},$$

by Theorem 1. Letting $\epsilon \to 0$ and $k \to \infty$,

(22)
$$\limsup_{a \to \infty} P[r \leq R_a \leq s, \xi_{t_a} \in I, \zeta_{t_a} \in J]$$
$$\leq (\frac{s-r}{\mu}) P[\inf_{j \leq -1}(Z_0^1 - Z_{-j}^1) \geq r, \xi_0 \in I] L\{J\}$$

By partitioning an interval $(r, s)$ into subintervals $r = r_0 < r_1 < \cdots < r_m = s$, applying (22) to each subinterval, and letting the partition become infinitely fine, it yields

$$\lim_{a \to \infty} P[r \leq R_a \leq s, \xi_{t_a} \in I, \zeta_{t_a} \in J] = \frac{1}{\mu} \int_r^s P[\inf_{j \leq -1} Z_j^* \geq u, \xi_0 \leq y] du \times L(z),$$

from which the theorem follows by using (21). □

Thus, the asymptotic distribution functions of $R_a$ and $\xi_{t_a}$ are

$$\frac{1}{\mu} \int_0^r P[\inf_{j \leq -1} Z_j^* \geq s] ds \qquad \text{and} \qquad \frac{1}{\mu} \int_{\xi_0 \leq y} [\inf_{j \leq -1} Z_j^*]_+ dP.$$

Denote the means of these distributions by $\rho$ and $\nu$.

**Theorem 4.** *If $p_x > 2/(1-q)$, then*

$$\lim_{a \to \infty} E|\frac{t_a}{a} - \frac{1}{\mu}| = 0. \tag{23}$$

*If also $p_y > (2-q)/q$ and (6) holds with $\kappa = 1$,*

$$E(t_a) = \frac{a + \rho - \nu - \lambda}{\mu} + o(1) \tag{24}$$

*as $a \to \infty$, where $\lambda = \lambda_1 + \cdots + \lambda_d$, and $\lambda_1, \ldots, \lambda_d$ are as in (8).*

*Proof.* From Lemma 1, $t_a/a \to^p 1/\mu$ as $a \to \infty$, and

$$\int_{t>M} t \, dP \leq M P[t > M] + \sum_{n=M}^{\infty} P[t > n] \to 0. \tag{25}$$

Relation (23) follows.

For (24), first observe that $P(A') = o(1/a)$, by Lemmas 1 and 3. From Wald's Lemmas and (23), $E(S_{t_a}) = \mu E(t_a)$ and $E[(S_{t_a} - \mu t_a)^2] = \sigma^2 E(t_a) = O(a)$. So, by (23),

$$\mu \int_A t_a dP = \int_A S_{t_a} dP + \int_{A'} (S_{t_a} - \mu t_a) dP,$$

and

$$|\int_{A'} (S_{t_a} - \mu t_a) dP| \leq \sqrt{P(A')} \sqrt{E(S_{t_a} - \mu t_a)^2} = o(1).$$

So, using Theorem 3 and Corollaries 1 and 2,

$$\mu \int_A t_a dP = aP(A) + \int_A (R_a - \xi_{t_a} - \zeta_{t_a}) dP + o(1) = a + \rho - \nu - \lambda + o(1).$$

Finally,

$$\int_{A'} t_a dP \leq MP(A') + \int_{t_a > M} t_a dP \to 0,$$

so that $E(t_a) = \int_A t_a dP + o(1)$. Relation (24) follows. $\square$

**Remark.** In fact, (23) can be established under the weaker condition $p_x > (2-q)/(1-q)$, at the expense of complicating the proof in fairly routine ways.

## 4. Examples

Consider the exponential model described in the Introduction. In this model, many examples are of the form

$$Z_n = ng(\frac{K_n}{n}, \frac{T_n^*}{n}), \tag{26}$$

where $g$ is a smooth function. This can be written in the form (1) with

$$S_n = ng(1, \frac{1}{\theta}) + g_{01}(1, \frac{1}{\theta}) \sum_{k=1}^n (L_k - \frac{1}{\theta}),$$

$$\xi_n = -\bigl[g_{10}(1,\tfrac{1}{\theta})\sum_{k=0}^{\infty}\mathbf{1}\{L_{n-k} > \sum_{j=0}^{k}\eta_{n-j}\} + g_{01}(1,\tfrac{1}{\theta})\sum_{k=0}^{\infty}\bigl(L_{n-k} - \sum_{j=0}^{k}\eta_{n-j}\bigr)_{+}\bigr]$$

and

$$\zeta_n = ng(\frac{K_n}{n}, \frac{T_n^*}{n}) - [S_n + \xi_n].$$

where $\eta_k = \tau_k - \tau_{k-1}$ and $g_{ij}$ denote the partial derivatives of $g$. The latter term is discussed in more detail below.

**Example 1** (Fixed width confidence intervals). Consider the problem of setting a confidence interval with fixed width $2h > 0$ for $\theta$ in the Exponential model. Asymptotic considerations lead to a stopping time of the form $t_a$ with

$$Z_n = n(\frac{T_n^*}{K_n})^2$$

and $a = c^2/h^2$, where $c$ is a normal percentile, as in [1]. This is of the form (26) with $g(x,y) = y^2/x^2$. □

**Example 2** (Repeated significance tests). The log-likelihood ratio test statistic for testing $\theta = 1$ is

$$Z_n = K_n \log\bigl(\frac{K_n}{T_n^*}\bigr) + (T_n^* - K_n),$$

which is of the form (26) with $g(x,y) = x\log(x/y) + (y-x)$. □

Returning to Equation (26), $\zeta_n$ may be written

$$\zeta_n = \zeta_{1,n} + \zeta_{2,n} + \zeta_{3,n},$$

where

$$\zeta_{1,n} = \frac{1}{2}n\bigl[g_{02}(1,\tfrac{1}{\theta})(\frac{T_n^*}{n} - \tfrac{1}{\theta})^2 + 2g_{11}(1,\tfrac{1}{\theta})(\frac{T_n^*}{n} - \tfrac{1}{\theta})(\frac{K_n}{n} - 1) + g_{20}(1,\tfrac{1}{\theta})(\frac{K_n}{n} - 1)^2\bigr],$$

$$\zeta_{2,n} = g_{10}(1,\tfrac{1}{\theta})\sum_{k=n}^{\infty}\mathbf{1}\{L_{n-k} > \sum_{j=0}^{k}\eta_{n-j}\} + g_{01}(1,\tfrac{1}{\theta})\sum_{k=n}^{\infty}\bigl(L_{n-k} - \sum_{j=0}^{k}\eta_{n-j}\bigr)_{+},$$

and

$$\zeta_{3,n} = Z_n - (S_n + \xi_n + \zeta_{1,n} + \zeta_{2,n}),$$

which may be bounded by a constant multiple of $n[|T_n^*/n - 1/\theta|^3 + |K_n/n - 1|^3]$ if $|\theta T_n^*/n - 1| \leq 1/2$ and $K_n/n \geq 1/2$. Further, letting $Y_k = L_k - 1/\theta$ and $T_n = Y_1 + \cdots + Y_n$,

$$\frac{T_n^*}{n} - \frac{1}{\theta} = \frac{T_n}{n} - \frac{\xi_n^o}{n}$$

where

$$\xi_n^o = \sum_{k=1}^{n}[L_{n-k+1} - (\tau_n - \tau_{n-k})]_{+}.$$

So, $\zeta_n$ may be written in the form (4) with

$$\zeta_n' = \frac{1}{2n}g_{02}(1,\tfrac{1}{\theta})T_n^2$$

and

$$\zeta_n'' = \zeta_{1,n}'' + \zeta_{2,n} + \zeta_{3,n}$$

with

$$\zeta''_{1,n} = -g_{02}(1, \frac{1}{\theta})\frac{T_n \xi_n^o}{n} + \frac{1}{2}g_{02}(1, \frac{1}{\theta})\frac{(\xi_n^o)^2}{n}$$
$$+ ng_{11}\left(1, \frac{1}{\theta}\right)\left(\frac{T_n}{n} - \frac{\xi_n^0}{n}\right)\left(\frac{K_n}{n} - 1\right) + \frac{n}{2}g_{20}\left(1, \frac{1}{\theta}\right)\left(\frac{K_n}{n} - 1\right)^2.$$

There are several terms to be considered in $\zeta''_n$, but none causes much difficulty. Some selected details: First $P[\zeta_{2,n} \neq 0] \leq P[L_k > \sum_{j=0}^{k} \eta_j$ for some $k \geq n]$, which is decays exponentially as $n \to \infty$. For $\zeta''_{1,n}$, $E|\xi_n^o|^p < \infty$ (independent of $n$) for all $p > 0$, and $E|T_n|^{2p} = O(n^p)$ for all $p > 0$. So, $E|\zeta''_{1,n}|^p = O(n^{-\frac{1}{2}p})$ for all $p$. Finally, let $B_n$ be the event that $|\theta T_n/n - 1| \leq 1/2$ and $K_n \geq n/2$. Then $P(B'_n)$ decays geometrically, and if $B_n$ occurs, then $\zeta_{3,n}$ is bounded by a constant multiple of

$$\frac{1}{n^2}\left[|T_n^* - \frac{n}{\theta}|^3 + |K_n - n|^3\right].$$

The conditions (6) and (7) can be verified using these facts.

## 5. Appendix

The Baum–Katz Inequalities are used in the following form: If $\lambda > 0$ and $n \geq l \geq 1$ are positive integers, then

$$P[\max_{k \leq n} |S_k - k\mu| > \lambda] \leq n\bar{F}(\frac{\lambda}{2l}) + \left(\frac{4l^2 n\sigma^2}{\lambda^2}\right)^l,$$

where $\bar{F}(\lambda) = P[|X_1| \geq \lambda]$. See [2] (pp. 373–375).

First consider $\Delta_0$. For any fixed $n$, $P[Z_n > a] \to 0$ as $a \to \infty$. So, it suffices to consider $\sum_{n=\ell}^{m} P[Z_n > a]$ for any $\ell$. If $n \leq m$, then

$$P[Z_n > a] \leq P[S_n - n\mu > \frac{a - n\mu}{2}] + P[\xi_n + \zeta_n > \frac{a - n\mu}{2}].$$

By (7), $C = \sup_{n \geq 0} E|\xi_n + \zeta_n|^2 < \infty$. So,

$$\sum_{n=n_0}^{m} P[\xi_n + \zeta_n > \frac{a - n\mu}{2}] \leq \sum_{n=n_0}^{m} \frac{4C}{(a - n\mu)^2} = O(a^{q-1}) \to 0$$

as $a \to \infty$. Next, let $l > 1/(1 - 2q)$. Then

$$\sum_{n=n_0}^{m} P[S_n - n\mu > \frac{a - n\mu}{2}] \leq \sum_{n=1}^{m} n\bar{F}(\frac{a - n\mu}{2l}) + m\left(\frac{16\ell^2 a\sigma^2}{a^{2(1-q)}}\right)^l,$$

which approaches 0 as $a \to \infty$, after some simple analysis. The analysis of $\Delta_1$ is similar, using

$$P[Z_n \leq a + \frac{1}{2}a^{1-q}] \leq P[|S_n - n\mu| > \frac{n\mu - a}{4}] + P[|\xi_n + \zeta_n| > \frac{n\mu - a}{4}]$$
$$\leq n\bar{F}(\frac{n\mu - a}{4}) + \left[\frac{16\ell^2 M\sigma^2}{(n\mu - a)^2}\right]^\ell + \frac{16C}{(n\mu - a)^2}.$$

and $P[t_a > n] \leq P[Z_n \leq a]$ for $n \geq M$. Relation (13) follows.

For (14),

$$aP[t_a > M] \leq aM\bar{F}\Big(\frac{a^{1-q}}{8}\Big) + a\Big(\frac{32\ell^2 M\sigma^2}{a^{2(1-q)}}\Big)^\ell + \frac{32aC}{a^{2(1-q)}},$$

which approaches 0 as $a \to \infty$ if $p_x > 2/(1-q)$ and $\ell > 1/(1-2q)$. That leaves $P[t_a \leq m]$. By (7), $C' = \sup_{n \geq n_0} E|\xi_n + \zeta_n|^\gamma < \infty$ for some $\gamma > (2-q)/(1-q)$. So,

$$P[t_a \leq m] \leq P[\max_{k \leq m}|S_k - k\mu| > \frac{1}{2}a^{1-q}] + \sum_{n=n_0}^{m} P[|\xi_n + \zeta_n| > \frac{1}{2}(a - n\mu)]$$

$$\leq m\bar{F}\Big(\frac{a^{1-q}}{4l}\Big) + \Big(\frac{4^{1+\gamma}l^2 m\sigma^2}{a^{2(1-q)}}\Big)^l + \sum_{n=n_0}^{m} \frac{C'}{(a-m\mu)^\gamma},$$

which is $o(1/a)$ if $l > 1/(1-2q)$. □